\documentclass[graybox]{svmult}

\usepackage{amsmath, amssymb}
\usepackage{amsfonts}
\usepackage{mathrsfs}
\usepackage[arrow,matrix,curve,cmtip,ps]{xy}

\usepackage{mathptmx}       
\usepackage{helvet}         
\usepackage{courier}        
\usepackage{type1cm}        
\usepackage{graphicx}        
\usepackage[bottom]{footmisc}

\newcommand{\coker}{\operatorname{coker}}         
\newcommand{\sgn}{\operatorname{sgn}}    
\newcommand{\Z}{\operatorname{\mathbb{Z}}}

\begin{document}

\title*{Classification of graph algebras: A selective survey}
\author{Mark Tomforde}
\institute{ Mark Tomforde, 
Department of Mathematics, University of Houston, Houston, TX 77204, USA \email{tomforde@math.uh.edu} \\
The final publication is available at link.springer.com. DOI forthcoming.
}
%
%
\maketitle

\abstract{This survey reports on current progress of programs to classify graph $C^*$-algebras and Leavitt path algebras up to Morita equivalence using $K$-theory.  Beginning with an overview and some history, we trace the development of the classification of simple and nonsimple graph $C^*$-algebras and state theorems summarizing  the current status of these efforts.  We then discuss the much more nascent efforts to classify Leavitt path algebras, and we describe the current status of these efforts as well as outline current impediments that must be solved for this classification program to progress.  In particular, we give two specific open problems that must be addressed in order to identify the correct $K$-theoretic invariant for classification of simple Leavitt path algebras, and we discuss the significance of various possible outcomes to these open problems.
}

\section{Introduction}

In 1976 Elliott proved a result, now known as Elliott's theorem, which states that direct limits of semisimple finite-dimensional algebras may be classified up to isomorphism by the dimension group (later identified with the scaled, ordered $K_0$-group) of the algebra.  Elliott's theorem implies that for the class of AF-algebras (i.e., $C^*$-algebraic direct limits of finite-dimensional $C^*$-algebras) the scaled, ordered topological $K_0$-group is a complete isomorphism invariant, and it also implies that in the class of ultramatricial algebras (i.e., algebraic direct limits of direct sums of semisimple finite-dimensional algebras over a fixed field) the scaled, ordered algebraic $K_0$-group is a complete isomorphism invariant.  (In this case the algebraic $K_0$-group and topological $K_0$-group coincide.)  Based on this result, Elliott boldly proposed that many additional classes of ($C^*$-)algebras can be classified up to isomorphism by $K$-theoretic invariants, and he famously formulated what is now called Elliott's conjecture: ``All separable nuclear simple $C^*$-algebras are classified up to isomorphism by $K$-theoretic invariants."  Work on this conjecture has been referred to as the Elliott program, and over the past four decades there have been numerous contributions made by several mathematicians, not the least of which are due to Elliott himself.  Very recently, work of Tikuisis, White, and Winter has completed the final steps for classifying all unital, separable, simple, and nuclear $C^*$-algebras of finite nuclear dimension which satisfy the UCT \cite{TWW}.  

The established work on the Elliott program is vast --- indeed, papers on the subject total thousands of pages and even a survey of all the accomplishments over the past four decades would most likely require a document the size of a book if it wished to be comprehensive.  (See the book \cite{Ror-book} for an introductory survey of the Elliott program with an emphasis on providing a technical overview of the Kirchberg-Phillips classification theorem, and see the papers \cite{ET, Ror} for a summary of the accomplishments in the Elliott program over the past 15 years.)

In the discussions here we wish to focus on two natural extensions of the Elliott program:

\begin{description}
\item[\textbf{Extension~1:}] Go beyond the simple $C^*$-algebras and attempt to classify certain nonsimple nuclear $C^*$-algebras using $K$-theory.
\smallskip
\smallskip
\item[\textbf{Extension~2:}] Step outside of the class of $C^*$-algebras and attempt to classify certain simple algebras using (algebraic) $K$-theory.  
\end{description}

\noindent Of course there is little to no hope these two extensions can be accomplished for all nuclear $C^*$-algebras or all simple algebras, so one major component of each program is to identify classes of nonsimple $C^*$-algebras and simple algebras that are amenable to classification.  Another major component of each program is determining exactly what $K$-theoretic data is needed for the classifying invariant.  

Significant progress for Extension~1 has been made for the class of graph $C^*$-algebras, while progress for Extension~2 has been more difficult, but had some stunning successes for the class of Leavitt path algebras.  These two classes, which we collectively refer to as \emph{graph algebras}, will be the focus of this survey.  

Readers who have no prior experience with graph algebras may initially (and incorrectly!) believe these classes are fairly small and specialized.  However, keep in mind that every AF-algebra is Morita equivalent to a graph $C^*$-algebra \cite{Dri} and every ultramatricial algebra over $\mathbb{C}$ is Morita equivalent to a Leavitt path algebra.  Thus the graph $C^*$-algebras and Leavitt path algebras are generalizations of the AF-algebras and ultramatrical algebras to which Elliott's theorem from 1976 applies.  As such, they are very suitable classes to explore at the beginnings of these programs.

Moreover, every Kirchberg algebra with free $K_1$-group is Morita equivalent to a graph $C^*$-algebra \cite{Szy}, so the class of graph $C^*$-algebras also generalizes many of the simple $C^*$-algebras to which the Kirchberg-Phillips classification theorem applies \cite{Phi}.  Consequently, the graph $C^*$-algebras comprise a large class containing several nonsimple $C^*$-algebras as well as many simple $C^*$-algebras of both AF and purely infinite type.  At the same time, the graph $C^*$algebras are not too large to escape classification.  The proposed invariant, called the filtered (or sometimes ``filtrated") $K$-theory, is a natural generalization of the invariant used for simple $C^*$-algebras, and contains the collection of all ordered $K_0$-groups and $K_1$-groups of subquotients of the $C^*$-algebra taking into account all the natural transformations among them.  While this seems like the obvious choice for the invariant, Meyer and Nest have constructed two separable, purely infinite $C^*$-algebras in the bootstrap class (each with a primitive ideal space having four points) that have the same filtered $K$-theory but are not Morita equivalent, thus demonstrating that filtered $K$-theory is inadequate to classify general nonsimple nuclear $C^*$-algebras.  This example also means that restricting to the class of graph algebras is not an artificial choice, but done out of necessity to obtain working theorems.  The counterexamples of Meyer and Nest lie outside the class of graph $C^*$-algebras, and the current working conjecture is that filtered $K$-theory suffices to classify graph $C^*$-algebras.

While a promising candidate for the classifying invariant of graph $C^*$-algebras has been identified, and many initial cases have been established successfully (see \cite{ERR} for a taxonomy), classification results for Leavitt path algebras have been more piecemeal, and the correct invariant for classification is still uncertain.  In this survey, we give a summary of the current status of classification results for simple and nonsimple graph $C^*$-algebras, and we describe how these results have guided initial work on classifying the simple Leavitt path algebras.  We also outline the existing classification results for simple Leavitt path algebras, and describe the current search for a complete Morita equivalence invariant in terms of $K$-theory.  We state two important open questions currently facing the classification program for simple Leavitt path algebras, and we also discuss the implications of various answers to these two open questions.

The ``selective survey" of the title refers to the fact that our attention will be primarily be concentrated on the ``geometric" classifications described in terms of moves on the graphs.  We will discuss the role that dynamical systems have played in these geometric classifications, and outline how dynamics results have been applied to graph algebras.  We will omit proofs in favor of focusing on the big picture, but do our best to explain the key ideas used to obtain results.  In addition, to avoid getting bogged down in too many technicalities, throughout this survey we shall restrict our attention to classification up to Morita equivalence (eschewing any mention of results for classification up to isomorphism).  

In addition to giving an update on current research, the author believes that the narrative, like many instances of mathematical investigation, provides an interesting story of the twists and turns that have occurred as several researchers have contributed to an area of investigation.

\section{Preliminaries}

In this section we establish notation and state some standard definitions.  To begin, we mention that we shall allow infinite graphs, but work under the standing hypothesis that all our graphs are countable.

\begin{definition}
A \emph{graph} $(E^0, E^1, r, s)$ consists of a countable set $E^0$ of vertices, a countable set $E^1$ of edges, and maps $r : E^1 \to E^0$ and $s : E^1 \to E^0$ identifying the range and source of each edge.  A graph is \emph{finite} if both the vertex set $E^0$ and the edge set $E^1$ are finite.
\end{definition}

\noindent Let $E := (E^0, E^1, r, s)$ be a graph. We say that a vertex $v
\in E^0$ is a \emph{sink} if $s^{-1}(v) = \emptyset$, and we say
that a vertex $v \in E^0$ is an \emph{infinite emitter} if
$|s^{-1}(v)| = \infty$.  A \emph{singular vertex} is a vertex that
is either a sink or an infinite emitter, and we denote the set of
singular vertices by $E^0_\textnormal{sing}$.  We also let
$E^0_\textnormal{reg} := E^0 \setminus E^0_\textnormal{sing}$, and
refer to the elements of $E^0_\textnormal{reg}$ as \emph{regular
vertices}; i.e., a vertex $v \in E^0$ is a regular vertex if and
only if $0 < |s^{-1}(v)| < \infty$.   

\subsection{Definitions of Graph Algebras}

Graph $C^*$-algebras were introduced in the 1990s, motivated by (and generalizing) earlier constructions, such as the Cuntz algebras and the Cuntz-Krieger algebras.  

\begin{definition}[The Graph $C^*$-algebra] \label{graph-C*-def}
If $E$ is a graph, the \emph{graph $C^*$-algebra} $C^*(E)$ is the universal
$C^*$-algebra generated by mutually orthogonal projections $\{ p_v
: v \in E^0 \}$ and partial isometries with mutually orthogonal
ranges $\{ s_e : e \in E^1 \}$ satisfying
\begin{enumerate}
\item $s_e^* s_e = p_{r(e)}$ \quad  for all $e \in E^1$
\item $s_es_e^* \leq p_{s(e)}$ \quad for all $e \in E^1$
\item $p_v = \sum_{\{ e \in E^1 : s(e) = v \}} s_es_e^* $ \quad for all $v \in E^0_\textnormal{reg}$.
\end{enumerate}
\end{definition}

\noindent Based on the success of graph $C^*$-algebras, in 2005 algebraists were inspired to define algebraic analogues, which are called Leavitt path algebras. 

\begin{definition}[The Leavitt Path Algebra] \label{LPA-lin-invo-def}
Let $E$ be a graph, and let $K$ be a field. We let $(E^1)^*$
denote the set of formal symbols $\{ e^* : e \in E^1 \}$.  The \emph{Leavitt path
algebra of $E$ with coefficients in $K$}, denoted $L_K(E)$,  is
the free associative $K$-algebra generated by a set $\{v : v \in
E^0 \}$ of pairwise orthogonal idempotents, together with a set
$\{e, e^* : e \in E^1\}$ of elements, modulo the ideal generated
by the following relations:
\begin{enumerate}
\item $s(e)e = er(e) =e$ for all $e \in E^1$
\item $r(e)e^* = e^* s(e) = e^*$ for all $e \in E^1$
\item $e^*f = \delta_{e,f} \, r(e)$ for all $e, f \in E^1$
\item $v = \displaystyle \sum_{\{e \in E^1 : s(e) = v \}} ee^*$ whenever $v \in E^0_\textnormal{reg}$.
\end{enumerate}
\end{definition}

\noindent As with the graph $C^*$-algebras, the Leavitt path algebras include many well-known classes of algebras and have been studied intensely in the algebra community since their introduction.  The interplay between these two classes of ``graph algebras" has been extensive and mutually beneficial --- graph $C^*$-algebra results have helped to guide the development of Leavitt path algebras by suggesting what results are true and in what direction investigations should be focused, and Leavitt path algebras have given a better understanding of graph $C^*$-algebras by helping to identify those aspects of $C^*(E)$ that are algebraic, rather than $C^*$-algebraic, in nature. Moreover, results from each class have had nontrivial applications to the other, and the work of researchers from each side has guided discovery for the other.  Indeed, nearly every theorem from each class seems to have a corresponding theorem in the other. For example, the graph-theoretic conditions on $E$ for which $C^*(E)$ is a simple algebra (respectively, an AF-algebra, a purely infinite simple algebra, an exchange ring, a finite-dimensional algebra) in the category of $C^*$-algebras are precisely
the same graph-theoretic conditions on $E$ for which $L_K(E)$ is a simple algebra (respectively, an ultramatricial
algebra, a purely infinite simple algebra, an exchange ring, a finite-dimensional algebra) in the category of $K$-algebras.  The exact reason for these similar properties is a bit of a mystery --- the graph $C^*$-algebra and Leavitt path algebra theorems are proven using different techniques, and the theorems for one class do not seem to imply the theorems for the other in any obvious way.  It has been suggested that there may exist some kind of ``Rosetta Stone" that would allow for translating and deducing one set of theorems from the other, but currently such a Rosetta Stone remains elusive.

We mention that when the underlying field of the Leavitt path is the complex numbers, then $L_\mathbb{C}(E)$ is isomorphic to a dense $*$-subalgebra of $C^*(E)$.  However, this alone is not enough to account for the similar results, since it is possible for dense $*$-subalgebras to have considerably different properties and structure from the ambient $C^*$-algebra.

It is also noteworthy that the field plays little role in most theorems for the Leavitt path algebras, and properties of $L_K(E)$ are frequently obtained entirely in terms of the graph $E$ with no dependence on the field $K$.  We will see in the final section of this survey that classification by $K$-theory is one of the few examples where this is not true, and the underlying field will be important in our theorems and in the invariants used for classification of simple Leavitt path algebras.

\subsection{Computation of $K$-groups}

We shall use the notation $K_n (A)$ for the $n$\textsuperscript{th} topological $K$-group of a $C^*$-algebra $A$, and we shall use the notation $K_n^\textnormal{alg} (R)$ for the $n$\textsuperscript{th} algebraic $K$-group of a ring $R$.  Due to a phenomenon called Bott periodicity, for any $C^*$-algebra $A$ we have $K_n (A) \cong K_{n+2} (A)$ for all $n \in \mathbb{Z}$.  Thus for $C^*$-algebras, all $K$-group information is contained in the $K_0$-group and $K_1$-group, and these are typically the only $K$-groups mentioned.  For rings (and Leavitt path algebras) there is no periodicity and all the algebraic $K$-groups may be distinct.

Computation of the topological $K$-groups of a graph $C^*$-algebra and the algebraic $K$-groups of a Leavitt path algebra is described in the following definition and proposition.

\begin{definition}
If $E = (E^0, E^1, r, s)$ is a graph, we define the vertex matrix $A_E$ to be the (possible infinite) square matrix indexed by the vertex set $E^0$, and for each $v, w \in E^0$ the entry $A_E(v,w)$ is equal to the number of edges in $E$ from $v$ to $w$.  Note that the entries of $A_E$ take values in $\{ 0, 1, 2, \ldots, \infty \}$.
\end{definition}

\noindent In the following proposition, for a set $X$ and an abelian group $G$, we use the notation $G^X := \bigoplus_{x \in X} G$.

\begin{proposition}[Computation of $K$-groups for Graph Algebras] \label{K-theory-comp-prop}
Let $E$ be a graph and decompose the vertices of $E$ as $E^0 = E^0_\textnormal{reg} \sqcup E^0_\textnormal{sing}$, and with respect to this decomposition write the vertex matrix of $E$ as $$A_E = \begin{pmatrix} B_E & C_E \\ * & * \end{pmatrix}$$ where $B_E$ and $C_E$ have entries in $\mathbb{Z}$ and the $*$'s have entries in $\mathbb{Z} \cup \{ \infty \}$.  For each $v \in E^0$, let $\delta_v \in \Z^{E^0}$ denote the vector with $1$ in the $v$\textsuperscript{th} position and $0$'s elsewhere, and for $x \in \Z^{E^0}$ let $[x]$ denote the equivalence class of $x$ in $\coker \left( \begin{pmatrix} I-B_E^t \\ -C_E^t\end{pmatrix} : \Z^{E^0_\textnormal{reg}} \to \Z^{E^0} \right)$.
\begin{itemize}
\item[(a)]  The topological $K$-theory of the graph $C^*$-algebra may be calculated as follows:  We have
$$K_0 (C^*(E)) \cong \coker \left( \begin{pmatrix} I-B_E^t \\ -C_E^t\end{pmatrix} : \Z^{E^0_\textnormal{reg}} \to \Z^{E^0} \right)$$ via an isomorphism that takes $[p_v]_0 \mapsto [\delta_v]$ and takes the positive cone of  $K_0 (C^*(E))$ to the cone of $\coker \left( \begin{pmatrix} I-B_E^t \\ -C_E^t\end{pmatrix} : \Z^{E^0_\textnormal{reg}} \to \Z^{E^0} \right)$ generated by the elements
$$
\{ [ p_v - \sum_{e \in F} s_es_e^*] : v \in E^0, F \subseteq s^{-1}(v), \text{ and $F$ finite} \},$$ and we have $$K_1 (C^*(E)) \cong \ker \left( \begin{pmatrix} I-B_E^t \\ -C_E^t\end{pmatrix} : \Z^{E^0_\textnormal{reg}} \to \Z^{E^0} \right).$$ 

\item[(b)] If $K$ is any field, then the algebraic $K$-theory of the Leavitt path algebra $L_K(E)$ may be calculated as follows:  We have
$$K_0^\textnormal{alg} (L_K(E)) \cong \coker \left( \begin{pmatrix} I-B_E^t \\ -C_E^t\end{pmatrix} : \Z^{E^0_\textnormal{reg}} \to \Z^{E^0} \right)$$ via an isomorphism that takes $[v]_0 \mapsto [\delta_v]$ and takes the positive cone of  $K_0^\textnormal{alg} (L_{K}(E))$ to the cone of $\coker \left( \begin{pmatrix} I-B_E^t \\ -C_E^t\end{pmatrix} : \Z^{E^0_\textnormal{reg}} \to \Z^{E^0} \right)$ generated by the elements
$$
\{ [ v - \sum_{e \in F} e e^*] : v \in E^0, F \subseteq s^{-1}(v), \text{ and $F$ finite} \},$$ and we have 
\begin{align*}
K_1^\textnormal{alg} (L_K(E)) \cong \ker \bigg( &\begin{pmatrix} I-B_E^t \\ -C_E^t\end{pmatrix} : \Z^{E^0_\textnormal{reg}} \to \Z^{E^0} \bigg) \\
& \oplus  \coker \left( \begin{pmatrix} I-B_E^t \\ -C_E^t\end{pmatrix} : (K_1^\textnormal{alg}(K))^{E^0_\textnormal{reg}} \to (K_1^\textnormal{alg} (K))^{E^0} \right)
\end{align*}
with $(K_1^\textnormal{alg}(K), +) \cong (K^\times, \cdot)$.
Moreover, there is a long exact sequence
$$\xymatrix{ K_n^\textnormal{alg}(K)^{E^0_\textnormal{reg}} \ar[r]^{ \left( \begin{smallmatrix} I-B_E^t \\ -C_E^t\end{smallmatrix} \right)} & K_n^\textnormal{alg}(K)^{E^0} \ar[r] & K_n^\textnormal{alg}(L_K(E)) \ar[r] & K_{n-1}^\textnormal{alg}(K)^{E^0_\textnormal{reg}}}
$$
for $n \in \Z$.
\end{itemize}
\end{proposition}

\noindent Remark: When $E$ has no singular vertices (which occurs, for example, whenever $E$ is finite with no sinks), then one may make the substitution $\left( \begin{smallmatrix} I-B_E^t \\ -C_E^t\end{smallmatrix} \right) = I - A_E^t$ in all the above expressions.

Throughout this survey we shall be concerned with classification of $C^*$-algebras and algebras up to Morita equivalence.  If $A$ is an algebra (or $C^*$-algebra) from a given class, and an object $\operatorname{I}(A)$ is assigned to $A$, we call the assignment a \emph{Morita equivalence invariant} for the class if 
$$
\text{$A$ Morita equivalent to $B$ $\implies$ $\operatorname{I}(A) = \operatorname{I}(B)$ \quad for all $A,B$ in the class}
$$
and 
we call the assignment a \emph{complete Morita equivalence invariant} for the class if 
$$
\text{$A$ Morita equivalent to $B$ $\iff$ $\operatorname{I}(A) = \operatorname{I}(B)$ \quad for all $A,B$ in the class.}
$$

\section{Classification of Simple and Nonsimple Graph $C^*$-algebras}

All graph $C^*$-algebras are nuclear and in the bootstrap class to which the UCT applies.  Furthermore, the  standing assumption that our graphs are countable ensures the associated graph $C^*$-algebras are separable.  Simple graph $C^*$-algebras are either AF (and classified by Elliott's theorem) or purely infinite (and classified by the Kirchberg-Phillips classification theorem).  Consequently, when $E$ is a graph and $C^*(E)$ is simple, the pair 
$$((K_0(C^*(E)), K_0(C^*(E))^+), K_1(C^*(E)))$$ 
consisting of the (pre-)ordered $K_0$-group together with the $K_1$-group is a complete Morita equivalence invariant for $C^*(E)$.  By looking at the ordering on the $K_0$-group, we can tell whether $C^*(E)$ is purely infinite or AF: If $K_0(C^*(E))^+ = K_0(C^*(E))$, then $C^*(E)$ is purely infinite and the invariant reduces to the pair $(K_0(C^*(E)), K_1(C^*(E))$; while if $K_0(C^*(E))^+ \subsetneq K_0(C^*(E))$, then $C^*(E)$ is AF, $K_1(C^*(E)) = 0$, and the ordered $K_0$-group $(K_0(C^*(E)), K_0(C^*(E))^+)$ classifies $C^*(E)$ up to Morita equivalence.

While the classification of simple graph $C^*$-algebras is a special case of the existing classification theorems for simple nuclear $C^*$-algebras (specifically Elliott's theorem and the Kirchberg-Phillips classification theorem), rather than merely considering the graph $C^*$-algebra result as a corollary, it perhaps better to take a historical view and think of simple graph $C^*$-algebras as the test cases that provided intermediate steps successively leading to more general theorems for classifying simple nuclear $C^*$-algebras.  Indeed, the classification of AF-algebras up to Morita equivalence given by Elliott's theorem, which initiated the entire classification program, is tantamount to classifying $C^*$-algebras of graphs with no cycles (since the two classes coincide up to Morita equivalence).  Likewise, the purely infinite simple graph $C^*$-algebras (especially particular subclasses) provided important initial steps leading to the Kirchberg-Phillips Classification theorem.  

For example, Cuntz first calculated the $K$-theory of the Cuntz algebra $\mathcal{O}_n$ (which is the $C^*$-algebra of a graph with one vertex and $n$ edges) in the paper \cite{Cun2}, showing that $K_0(\mathcal{O}_n) \cong \mathbb{Z} / n \mathbb{Z}$.  This implies a very specific case of the Kirchberg-Phillips theorem for the Cuntz algebras: $\mathcal{O}_m$ is Morita equivalent to $\mathcal{O}_n$ $\iff$ $K_0(\mathcal{O}_m) \cong K_0(\mathcal{O}_n)$ $\iff$ $m=n$.  

Likewise, the classification of simple Cuntz-Krieger algebras was an important early step in the classification program.  (The Cuntz-Krieger algebras are precisely the $C^*$-algebras of finite graphs with no sinks or sources.) Simple Cuntz-Krieger algebras are purely infinite and classified up to Morita equivalence by their $K_0$-group. (The $K_1$-group turns out to be redundant, because the $K_1$-group of a Cuntz-Krieger algebra is isomorphic to the free part of the $K_0$-group.)  The groundwork for this classification was laid by Cuntz and Krieger \cite{CK}, who recognized the connection with dynamics, and the final portions of the classification were later established by R\o rdam (using an important lemma outlined by Cuntz in a talk) \cite{Ro}.  We will return to this result, its proof, and the connection with dynamics in the next section.

In addition to providing stepping stones toward more general classification theorems for simple nuclear $C^*$-algebras, the graph $C^*$-algebras have also provided a class for exploring classification of nonsimple $C^*$-algebras.  In a graph $C^*$-algebra with finitely many ideals, each quotient and each ideal is Morita equivalent to a graph $C^*$-algebra, so any graph $C^*$-algebra with finitely many ideals may be built up from the simple graph $C^*$-algebras by taking extensions a finite number of times. 

The invariant used to classify a nonsimple graph $C^*$-algebra $C^*(E)$ is called the filtered (or sometimes ``filtrated") $K$-theory and denoted $\textrm{FK}_X^+ (C^*(E))$, where $X$ is the primitive ideal space of $C^*(E)$.  The filtered $K$-theory contains the collection of all ordered $K_0$-groups and $K_1$-groups of subquotients of the $C^*$-algebra, taking into account all the natural transformations among them (see \cite[\S1]{ERR} for a precise definition).  When $C^*(E)$ has a single (nonzero) ideal, the filtered $K$-theory is simply the six-term exact sequence in $K$-theory (including the ordering on all $K_0$-groups) determined by the unique ideal.  Eilers and the author \cite{ET1} proved that the filtered $K$-theory (i.e., the six-term exact sequence) is a complete Morita equivalence invariant for the class of graph $C^*$-algebras with one ideal.  This kicked off a flurry activity in which several researchers established classification results for graph $C^*$-algebras with multiple ideals.  A summary of these results and description of the status quo for this program can be found in the survey \cite{ERR} --- in particular, the authors there describe how we have a complete classification of graph $C^*$-algebras with a primitive ideal space having three or fewer points, and for graph $C^*$-algebras whose primitive ideal space has four points 103 of the 125 cases have been solved, leaving less than one fifth of the cases open.  In addition, as we will discuss at the end of Section~\ref{moves-sec}, Eilers, Restorff, and Ruiz recently announced that they have recently shown that filtered $K$-theory is a complete Morita equivalence invariant for all unital graph $C^*$-algebras.

\section{In Search of Techniques to Classify Algebras: Shift Spaces, Flow Equivalence, and Moves on Graphs}
\label{moves-sec}

Many of the techniques used to establish classification theorems for $C^*$-algebras have no hope of going through for general algebras.  Indeed, the classification proofs frequently make use of the $C^*$-algebra structure (e.g., the completeness is used frequently to take limits).  As one undertakes a classification of algebras, the first step is to begin with simple algebras.  In addition, it is sensible to restrict initial attention to algebras that are somehow similar to the  $C^*$-algebras for which the classification program had early successes using more modest methods.  The graph $C^*$-algebras (particularly the Cuntz-Krieger algebras) are such a class, and hence their algebraic analogues, the Leavitt path algebras (particularly Leavitt path algebras of finite graphs), arise as natural candidates for attempts at classification.

Simple Leavitt path algebras exhibit a dichotomy similar to simple graph $C^*$-algebras: A simple Leavitt path algebra is either ultramatricial (if the graph has no cycles) or purely infinite (if the graph has a cycle).  In the ultramatricial case, Elliott's theorem applies and the Leavitt path algebra is classified by its ordered $K_0$-group.  In the purely infinite case, we seek an algebraic analogue of the Kirchberg-Phillips classification theorem, and therefore we look at how the classification was obtained for early  investigations into special cases --- more specifically, we shall carefully examine the classification of simple Cuntz-Krieger algebras.  

The Cuntz-Krieger algebras correspond to $C^*$-algebras of finite graphs with no sinks and no sources, and the simplicity of the Cuntz-Krieger algebras corresponds to the graph being strongly connected (i.e., for each pair of vertices $v$ and $w$ there is a path from $v$ to $w$ and a path from $w$ to $v$) and not a single cycle.  Let $E = (E^0, E^1, r, s)$ be a finite graph with no sinks and no sources.  One may define the (two-sided) shift space $$X_E := \{ \ldots e_{-2}e_{-1}.e_0 e_1 e_2 \ldots  \ : \  e_i \in E^1 \text{ and } r(e_i) = s(e_{i+1}) \text{ for all } i \in \mathbb{Z} \}$$
consisting of all bi-infinite paths in the graph, together with the shift map $\sigma : X_E \to X_E$ given by $\sigma(  \ldots e_{-2}e_{-1}.e_0 e_1 e_2 \ldots) =  \ldots e_{-1}e_{0}.e_1 e_2 e_3 \ldots$.  Cuntz and Krieger observed a connection between the (characterizations of) flow equivalence of this dynamical system and the Morita equivalence class of the $C^*$-algebra associated with the graph.

If $E$ and $F$ are finite graphs with no sinks and no sources, the shift spaces $X_E$ and $X_F$ are \emph{flow equivalent} if their suspension flows are homeomorphic via a homeomorphism that carries orbits to orbits and preserves each orbit's orientation.  A precise definition of flow equivalence may be found in \cite[\S13.6]{LM}, but as we shall see shortly, for the purposes of this survey one does not need to understand flow equivalence so much as note that it is equivalent to other conditions.

Franks gave an algebraic characterization of flow equivalence for strongly connected graphs: If $E$ and $F$ are strongly connected finite graphs, then $X_E$ and $X_F$ are flow equivalent if and only if $\coker (1-A^t_E) = \coker (1-A^t_F)$ and $\sgn (\det (1-A^t_E)) = \sgn (\det (1-A^t_F)$.  Here $A_E$ is the vertex matrix of the graph $E$, and $\sgn (\det (1-A^t_E))$ is the sign of the number $\det (1-A^t_E)$; i.e., the value $+$, $-$, or $0$.

In addition, Parry and Sullivan gave a different characterization of flow equivalence based on ``moves"; i.e., operations that may be performed on the graph and which preserve flow equivalence of the associated shift space.  Parry and Sullivan needed three moves for their characterization, which are named as follows:

$ $

\textbf{Move (O):} Outsplitting \qquad \textbf{Move (I):} Insplitting \qquad \textbf{Move (R):} Reduction

$ $

Precise definitions of these moves can be found in \cite[\S3]{RT}, but for the purposes of this survey any readers unfamiliar with the moves may be better served by informal descriptions: \emph{Outsplitting} allows one to partition the outgoing edges of a vertex into nonempty sets and ``split" the vertex and incoming edges so that each set of outgoing edges from the partition now emits from its own vertex.  \emph{Insplitting} allows one to partition the ingoing edges of a vertex into nonempty sets and ``split" the vertex and outgoing edges so that each set of ingoing edges from the partition now enters its own vertex.   \emph{Reduction} allows one to ``collapse" certain vertices that have a single edge going from one vertex to the other.  

For each move there is also an \emph{inverse move}, so that if a Move $X$ is applied to the graph $E$ to obtain the graph $E'$, we say $E$ is obtained by performing the inverse of Move $X$ to $E'$.  Although we will not need their names, for the reader's edification we mention that the inverse of outsplitting is called \emph{outamalgamation}, the inverse of insplitting is called \emph{inamalgamation}, and the inverse of reduction is called \emph{delay}.  

Parry and Sullivan proved that if $E$ and $F$ are strongly connected graphs, then $X_E$ and $X_F$ are flow equivalent if and only if the graph $E$ may turned into the graph $F$ by finitely many applications of Moves (O), (I), (R), and their inverses.

Combining Franks result with the result of Parry and Sullivan, we thus obtain the following:

\begin{theorem}[Franks, Parry and Sullivan] \label{flow-equivalence-thm}
Let $E$ and $F$ be strongly connected finite graphs.  Then the following are equivalent:
\begin{itemize}
\item[(1)] \ The shift spaces $X_E$ and $X_F$ are flow equivalent.
\item[(2)] \ $\coker (1-A^t_E) = \coker (1-A^t_F)$ and $\sgn (\det (1-A^t_E)) = \sgn (\det (1-A^t_F))$.
\item[(3)] \ The graph $E$ may turned into the graph $F$ by finitely many applications of Moves (O), (I), (R), and their inverses.
\end{itemize}
\end{theorem}

Cuntz and Krieger had multiple insights to recognize the relationship of flow equivalence with the Morita equivalence of Cuntz-Krieger algebras.  First, after calculating the $K$-groups of a Cuntz-Krieger algebra, Cuntz and Krieger observed that the $K_0$-group coincides with the group $\coker (1-A^t_E)$ appearing in the flow equivalence classification.  (This group is sometimes called the Bowen-Franks group by dynamicists.)  The second important observation of Cuntz and Krieger was that the Moves (O), (I), (R) (and consequently their inverses) preserve Morita equivalence of the associated $C^*$-algebra.  

These observations, combined with Theorem~\ref{flow-equivalence-thm}, imply that if $E$ and $F$ are strongly connected finite graphs with isomorphic $K_0$-groups and with $\sgn (\det (1-A^t_E)) = \sgn (\det (1-A^t_F))$, then $E$ can be transformed into $F$ using a finite number of Moves (O), (I), (R) and their inverses, and consequently the $C^*$-algebras of $E$ and $F$ are Morita equivalent.  (Note that we really only need the equivalence of (2) and (3) in Theorem~\ref{flow-equivalence-thm}, and for the purposes of Cuntz and Krieger's result, we can completely ignore the notion of flow equivalence if we wish, viewing (2) $\iff$ (3) as a purely combinatorial fact about graphs.)

Although Cuntz and Krieger formulated their study of the Cuntz-Krieger algebras in terms of matrices, we wish to use the more modern approach of describing the $C^*$-algebras in terms of graphs.  The Cuntz-Krieger algebras may be thought of as the $C^*$-algebras of finite graphs with no sinks or sources.  If a graph $C^*$-algebra is simple and purely infinite, then the graph cannot contain a sink, so when we restrict to the simple purely infinite case we automatically have the ``no sinks" condition for Cuntz-Krieger algebras.  However, we do need a method to deal with sources, and to accomplish this we introduce a new move:

\begin{center}
\textbf{Move (S):} Source Removal
\end{center}

\noindent A precise definition of Move (S) can be found in \cite[\S3]{RT}, but an informal description is fairly accurate and informative: To perform Move (S) we select a source vertex in the graph and then remove this vertex and all edges beginning at this vertex.  As with the other moves, Move (S) preserves Morita equivalence of the associated $C^*$-algebra.  The is also an inverse move called \emph{source addition}.  Note that the process of performing Move (S) removes a source, but may create other sources in doing so.  Nonetheless, one can easily show that in a finite graph with no sinks, repeated applications of Move (S) will ultimately (and in a finite number of steps) result in a graph with no sources.

The following is a reformulation of Cuntz and Krieger's result in the language of graphs that also takes the presence of sources into account.

\begin{theorem}[Cuntz and Krieger] \label{sufficient-CK-thm}
Let $E$ and $F$ be finite graphs for which $C^*(E)$ and $C^*(F)$ are simple and purely infinite.  If $K_0(C^*(E)) \cong K_0(C^*(F))$ and $\sgn (\det (1-A^t_E)) = \sgn (\det (1-A^t_F))$, then $C^*(E)$ is Morita equivalent to $C^*(F)$.  Moreover, in this case, the graph $E$ may turned into the graph $F$ by finitely many applications of Moves (S), (O), (I), (R), and their inverses.
\end{theorem}

Cuntz and Krieger suspected that $\sgn (\det (1-A^t_E))$ was not a necessary condition for Morita equivalence of the $C^*$-algebras, but they were unable to remove the hypothesis from their theorem.  It was not until 15 years later that R\o rdam was able to remove the ``sign of the determinant" condition and obtain a complete Morita equivalence invariant for Cuntz-Krieger algebras.  To accomplish this, R\o rdam used an additional graph move that did not appear in the study of flow equivalence.  This move is called the \emph{Cuntz splice}, and because it will be important for us in the remainder of this survey, we shall describe it in greater detail than the other moves.

$$
\text{\textbf{Move (CS):} Cuntz Splice}
$$

\noindent  If $E$ is a graph and $v$ is any vertex in $E$ that is the base of two distinct cycles, then Move (CS) is performed by "splicing" on the following additional portion to $E$:

$ $

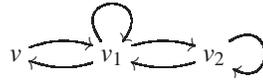
$$
\xymatrix{ 
		 v  \ar@/^0.5em/[r] & v_1 \ar@/^0.5em/[l] \ar@/^0.5em/[r] \ar@(ul,ur) & v_2 \ar@/^0.5em/[l] \ar@(ur,dr)
	}
$$

$ $

\noindent Here is an example showing the Cuntz splice performed on a graph with two vertices and three edges to produce a new graph with four vertices and nine edges.
 
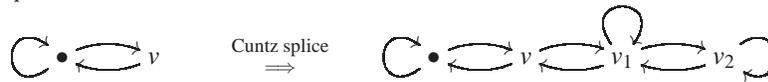
\begin{example}
$$
\xymatrix{ 
		\bullet \ar@(dl,ul) \ar@/^0.5em/[r] & v \ar@/^0.5em/[l]
	} 
\qquad \xymatrix{ {\text{Cuntz splice} \atop \Longrightarrow} }
\qquad \quad 
\xymatrix{ 
		\bullet \ar@(dl,ul) \ar@/^0.5em/[r] & v \ar@/^0.5em/[l] \ar@/^0.5em/[r] & v_1 \ar@/^0.5em/[l] \ar@/^0.5em/[r] \ar@(ul,ur) & v_2 \ar@/^0.5em/[l] \ar@(ur,dr)
	}
$$
\end{example}

$ $

\noindent Unlike the other moves, we shall have no need of an inverse move for the Cuntz splice.  The usefulness of the Cuntz splice is due to the following fact:  If $E$ is a graph and $E^-$ is obtained by performing a Cuntz splice to a vertex of $E$, then $K_0(C^*(E)) \cong K_0(C^*(E^-))$ and $\sgn (\det (I - A_E^t)) = -\sgn (\det (I - A_{E^-}^t))$.  In other words, the Cuntz splice preserves the $K_0$-group of the associated $C^*$-algebra while ``flipping" the sign of the determinant.  R\o rdam's main contribution was to prove that the Cuntz splice preserves Morita equivalence of the associated $C^*$-algebra.  Using this fact, one can start with two purely infinite simple graph $C^*$-algebras having the same $K_0$-group.  If the signs of the determinants are the same, simply apply Theorem~\ref{sufficient-CK-thm} to deduce Morita equivalence.  If not, apply the Cuntz splice once to one of the graphs to switch the sign of the determinant, and then apply Theorem~\ref{sufficient-CK-thm} to deduce Morita equivalence.  R\o rdam phrased his result in terms of the matrix description of Cuntz-Krieger algebras, but we reformulate it here in the modern  language of graphs, which also takes the presence of sources into account.

\begin{theorem}[R\o rdam] \label{ME-finite-thm}
Let $E$ and $F$ be finite graphs for which $C^*(E)$ and $C^*(F)$ are simple and purely infinite.  Then the following are equivalent:
\begin{itemize}
\item[(1)] \ $C^*(E)$ is Morita equivalent to $C^*(F)$.
\item[(2)] \ $K_0(C^*(E)) \cong K_0(C^*(F))$.
\item[(3)] \ The graph $E$ may turned into the graph $F$ by finitely many applications of Moves (S) (O), (I), (R), and their inverses, and at most one application of Move~(CS).  Moreover, no applications of Move~(CS) are needed if $\sgn (\det (1-A^t_E)) = \sgn (\det (1-A^t_F))$, and exactly one application of Move (CS) is required otherwise.
\end{itemize}
\end{theorem}

Theorem~\ref{ME-finite-thm} shows that $K_0(C^*(E))$ is a complete invariant for Morita equivalence in the class of simple purely infinite $C^*$-algebras of finite graphs.  Even better than that, Theorem~\ref{ME-finite-thm} gives moves on the graph generating the equivalence relation (in analogy to moves for other equivalence relations, such as the Reidemeister moves for the isotopy class of a knot).  This allows one to turn the question of Morita equivalence of the $C^*$-algebras into a combinatorial problem on graphs.  For this reason, the moves of Theorem~\ref{ME-finite-thm} are sometimes said to give a ``geometric classification" of these graph $C^*$-algebras.

In addition, for finite graphs whose $C^*$-algebras are purely infinite and simple, Theorem~\ref{flow-equivalence-thm} and Theorem~\ref{ME-finite-thm} explain the precise relationship between flow equivalence of the shift space and Morita equivalence of the $C^*$-algebra.  In particular, $K_0(C^*(E))$ is a complete invariant for Morita equivalence of $C^*(E)$, and the pair $(K_0(C^*(E)), \sgn (\det (1-A^t_E)))$ is a complete invariant for flow equivalence of $X_E$.  Consequently, for these graphs the flow equivalence of the shift space is a finer equivalence relation than Morita equivalence of the $C^*$-algebra (i.e.,  $X_E$ flow equivalent to $X_F$ implies that $C^*(E)$ is Morita equivalent to $C^*(F)$, but not conversely in general).

After R\o rdam's work and this geometric classification of Cuntz-Krieger algebras, efforts in the classification program progressed in ways that were less ``geometric".  However, in 2005 Abrams and Aranda Pino introduced Leavitt path algebras.  In 2008, approximately 13 years after the geometric classification for simple Cuntz-Krieger algebras was obtained, Abrams, Louly, Pardo, and Smith were inspired to seek a similar geometric classification for simple Leavitt path algebras of finite graphs \cite{ALPS}.  We will discuss the Leavitt path algebra classification in Section~\ref{LPA-sec}, and as with graph $C^*$-algebras we shall see the sign of the determinant condition is a stumbling block.  Unlike the $C^*$-algebra situation, however, this problem has not been resolved and it is an open question as to whether the sign of the determinant may be removed.  We will discuss the ramifications of this question, and the implications of possible answers, in the next section.

With Leavitt path algebra considerations causing attention to be returned to a geometric classification, S\o rensen had the novel idea to reconsider the graph $C^*$-algebras and seek geometric classifications for simple $C^*$-algebras of infinite graphs \cite{Sor}.  Although the Kirchberg-Phillips theorem, established in 2000, showed that the pair of the $K_0$-group and $K_1$-group is a complete Morita equivalence invariant for purely infinite simple graph $C^*$-algebras, the result is highly non-geometric and does not allow one to establish the Morita equivalence in any concrete or constructive way.  S\o rensen's key insight was to realize that classification by moves could still be obtained when the graph has a finite number of vertices and an infinite number of edges.  In this case, the $K$-groups of the $C^*$-algebras are obtained as the cokernel and kernel of a finite rectangular (but not square) matrix indexed by the vertices.  In this situation, due to the fact there are infinitely many edges, one cannot define a shift space as before.  This is because definition of a shift spaces requires a finite alphabet (i.e., finitely many edges) from which each position in the bi-infinite sequences may be chosen.  Despite this, one can still focus on the equivalence of (2) and (3) in Theorem~\ref{flow-equivalence-thm}, considering it as a purely combinatorial fact about graphs, and seek an analogous result for infinite graphs with finitely many vertices.  As one would expect, the object $\coker (I-A_E^t) \cong K_0(C^*(E))$ from the finite graph case must now be replaced by the pair $(K_0(C^*(E)), K_1(C^*(E))$.  It is less clear what quantity should play the analogous role of $\sgn (\det (I-A_E^t))$, since the matrix involved is not square and hence its determinant does not exist.  Surprisingly, S\o rensen proved that the sign of the determinant condition simply disappears in the presence of infinitely many edges, and --- in what is even better news --- this means there is no need for the Cuntz splice.  S\o rensen's results from \cite{Sor} may be summarized as follows.

\begin{theorem}[S\o rensen] \label{ME-infinite-thm}
Let $E$ and $F$ be graphs with a finite number of vertices and an infinite number of edges and with the property that $C^*(E)$ and $C^*(F)$ are simple.  Then the following are equivalent:
\begin{itemize}
\item[(1)] \ $C^*(E)$ is Morita equivalent to $C^*(F)$.
\item[(2)] \ $K_0(C^*(E)) \cong K_0(C^*(F))$ and $K_1(C^*(E)) \cong K_1(C^*(F))$.
\item[(3)] \ The graph $E$ may turned into the graph $F$ by finitely many applications of Moves (S), (O), (I), (R), their inverses.
\end{itemize}
\end{theorem}

The fact that S\o rensen's result does not involve the Cuntz splice has two important consequences: First, the moves (S), (O), (I), (R), and their inverses produce explicit Morita equivalences, and consequently one can concretely construct the imprimitivity bimodule linking $C^*(E)$ to $C^*(F)$ by using full corners of the $C^*$-algebras of the intermediary graphs between $E$ and $F$.  (The Cuntz splice does not produce an explicit Morita equivalence, so this concrete construction cannot be accomplished for the finite graphs in Theorem~\ref{ME-finite-thm} when $\sgn (\det (1-A^t_E)) \neq \sgn (\det (1-A^t_F))$.)  Second, when seeking a version of Theorem~\ref{ME-infinite-thm} for Leavitt path algebras, the sign of the determinant condition is no longer present to cause difficulties.  

A graph $C^*$-algebra $C^*(E)$ is unital precisely when the graph $E$ has a finite number of vertices.  Since any unital simple graph $C^*$-algebra is either Morita equivalent to $\mathbb{C}$ or purely infinite, the combination of Theorem~\ref{ME-finite-thm} and Theorem~\ref{ME-infinite-thm} gives a complete classification of unital simple graph $C^*$-algebras up to Morita equivalence.

In May 2015 Eilers, Restorff, Ruiz, and S\o rensen posted a preprint \cite{ERRS} to the arXiv in which they extended the geometric classification to all unital graph $C^*$-algebras of real rank zero.  Specifically, they show that filtered $K$-theory is a complete Morita equivalence invariant for unital graph $C^*$-algebras of real rank zero, and that when Morita equivalence occurs, one graph may be turned into the other using a finite number of Moves (S), (O), (I), (R), their inverses, and Move (CS).  In July 2015, Eilers, Restorff, Ruiz, and S\o rensen posted an update to their arXiv entry stating that they are now able to remove the hypothesis of real rank zero and give a geometric classification for all unital graph $C^*$-algebras.  They also stated that the preprint \cite{ERRS} will not be published, and instead a paper with the more general results (containing all results of \cite{ERRS} as special cases) will be written and published in its place (see \texttt{http://front.math.ucdavis.edu/1505.06773}).  In talks, Eilers has stated that the  filtered $K$-theory is a complete Morita equivalence invariant for unital graph $C^*$-algebras and that a geometric classification is possible.  However, when the unital graph $C^*$-algebra is not real rank zero, we must include one additional move besides Moves (S), (O), (I), (R), their inverses, and Move (CS), in order to handle the situation of cycles with no exits.  Since this additional move was discovered while Eilers and Restorff were visiting Ruiz at his home institution in Hawai'i, and since the move involves a graphical picture similar to a butterfly, the authors have tentatively called the move ``pulelehua" --- the Hawaiian word for butterfly.

\section{Classification of Leavitt Path Algebras of Finite Graphs}
\label{LPA-sec}

Since the Leavitt path algebras are defined in a manner analogous to the graph $C^*$-algebras, they are a natural candidate for a class of algebras that may be amenable to classification by $K$-theory.  In addition, the geometric classification of unital graph $C^*$-algebras, described in terms of moves on the graph, provides a viable approach to classification for Leavitt path algebras.  

Abrams, Louly, Pardo, and Smith initiated the classification of Leavitt path algebras in \cite{ALPS}.  (A preprint of \cite{ALPS} was posted to the arXiv in 2008, and a published version appeared in 2011.)  To begin, they observed that for any graph $E$ and any field $K$, one has $K_0^\textnormal{alg} (L_K(E)) \cong K_0(C^*(E))$ so that the algebraic $K_0$-group of the Leavitt path algebra agrees with the $K_0$-group of the graph $C^*$-algebra and is independent of the field $K$.  In addition, Abrams, Louly, Pardo, and Smith proved that the graph moves (S), (O), (I), (R), and their inverses preserve Morita equivalence of the Leavitt path algebra of the graph.  However, they were unable to determine whether or not the Cuntz splice preserves Morita equivalence of the associated Leavitt path algebra.  Thus, following the proof strategy established by Cuntz and Krieger (with later contributions by R\o rdam), they could not avoid the sign of the determinant condition and were only able to establish sufficient conditions for Morita equivalence.  We state their result here.

\begin{theorem}[Abrams, Louly, Pardo, and Smith]
Let $K$ be any field, and let $E$ and $F$ be finite graphs for which $L_K(E)$ and $L_K(F)$ are simple and purely infinite.  If $K_0^\textnormal{alg}(L_K(E)) \cong K_0^\textnormal{alg}(L_K(F))$ and $\sgn (\det (1-A^t_E)) = \sgn (\det (1-A^t_F)$, then $L_K(E)$ is Morita equivalent to $L_K(F)$.  Moreover, in this case, the graph $E$ may turned into the graph $F$ by finitely many applications of Moves (S), (O), (I), (R), and their inverses.
\end{theorem}

One noteworthy consequence of this result is that, as with most of the Leavitt path algebra results, the field plays no role in determining the Morita equivalence class of the Leavitt path algebra for these particular types of graphs.   Indeed, the invariant sufficient for classification, the pair $(K_0^\textnormal{alg}(L_K(E)), \sgn (\det (1-A^t_E)))$, depends only on the graph $E$ and is independent of the field.  

Currently, it is unknown whether the sign of the determinant is a Morita equivalence invariant for Leavitt path algebras.  This puts the classification program for simple Leavitt path algebras of finite graphs in a similar state that the classification of simple graph $C^*$-algebras found itself in during the 15 years period following Cuntz and Krieger's work and prior to R\o rdam's contribution of the Cuntz splice.  As a result, resolving whether the sign of the determinant can be removed is currently one of the central issues in the classification of simple Leavitt path algebras.  This problem is equivalent to determining whether the Cuntz splice preserves Morita equivalence, and thus the open question may be formulated as follows.

$ $

\noindent \textbf{Open Question 1:} \emph{Let $K$ be a field and let $E$ be a finite graph such $\det (1-A^t_E) \neq 0$ and $L_K(E)$ is simple and purely infinite.  If $E^-$ denotes the graph obtained by performing a Cuntz splice to $E$, then are the Leavitt path algebras $L_K(E)$ and $L_K(E^-)$ Morita equivalent?}

$ $

This question has been open since the first preprint of \cite{ALPS} appeared in 2008, and it is currently at the forefront of the classification program for Leavitt path algebras.  Many researchers have worked on this problem, but with little to show for their efforts.  In fact, we currently cannot answer the question in even elementary special cases.

For example, suppose we take $E_2$ to be the graph with one vertex and two edges (arguably, the most basic example of the graphs the question is asking about), and let $E_2^-$ be the graph obtained by performing a Cuntz splice at the vertex of $E$.

$ $

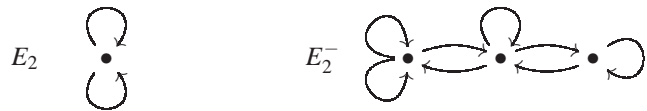
$$E_2 \qquad \xymatrix{ 
		\bullet \ar@(ul,ur) \ar@(dl,dr) 
	} 
	\qquad
	\qquad
	\quad
	\qquad E_2^- \qquad
	\xymatrix{ 
		\bullet \ar@(l,u) \ar@(l,d) \ar@/^0.5em/[r] & \bullet \ar@/^0.5em/[l] \ar@/^0.5em/[r] \ar@(ul,ur) & \bullet \ar@/^0.5em/[l] \ar@(ur,dr)
	}
$$

$ $

\noindent Then $L_K(E_2)$ is isomorphic to the Leavitt algebra $L_2$, and we have $K_0^\textnormal{alg} (L_K(E_2)) \cong K_0^\textnormal{alg} (L_K(E_2^-)) \cong \{ 0 \}$, $\det (1-A^t_{E_2}) = 1$, and $\det (1-A^t_{E_2^-})=-1$.  However, it is currently an open question as to whether $L_K(E_2)$ and $L_K(E_2^-)$ are Morita equivalent.

Besides restricting the graph, another approach to finding a more tractable special case of Open Question~1 is to restrict the field, for instance to fix $K = \mathbb{C}$ or $K = \mathbb{Z}_2$.  However, no results for this special case have been obtained either.  Even combinations of these restrictions (i.e., restricting both the graph and the field) yield open problems --- no one knows, for example, whether $L_\mathbb{C}(E_2)$ and $L_\mathbb{C}(E_2^-)$ are Morita equivalent, or whether $L_{\mathbb{Z}_2}(E_2)$ and $L_{\mathbb{Z}_2}(E_2^-)$ are Morita equivalent.

In July 2015, Johansen and S\o rensen announced the preprint \cite{JS}, which to the author's knowledge contain some of the first concrete results concerning the sign of the determinant condition.  Although Leavitt path algebras are defined over fields, as noted by the author in \cite{Tom}, for any graph $E$ and any commutative ring $R$ it is possible to construct a Leavitt path algebra $L_R(E)$ with coefficients in $R$.  Johansen and S\o rensen proved that if we choose the coefficients to be the ring $\mathbb{Z}$, then $L_{\mathbb{Z}}(E_2)$ is not $*$-isomorphic to $L_{\mathbb{Z}}(E_2^-)$.  (This is contrasted with the graph $C^*$-algebra situation, where $C^*(E_2)$ is $*$-isomorphic to $C^*(E_2^-)$.)  Consequently, the Cuntz splice does not preserve $*$-isomorphism of Leavitt path algebras over the ring $\mathbb{Z}$.  What this means --- if anything --- for Morita equivalence (instead of $*$-isomorphism) of Leavitt path algebras over fields (instead of rings), is yet unclear.  But at the very least, Johansen and S\o rensen's result shows us that not all Cuntz splice results for graph $C^*$-algebras will generalize to algebras over commutative rings, and this raises the potential for some unexpected phenomena with Leavitt path algebras over fields. More importantly, up to this point a preponderance of researchers' efforts have been spent trying to prove that the Cuntz splice does preserve Morita equivalence of Leavitt path algebras (over fields).  Johansen and S\o rensen's result suggests that perhaps we should be spending more time trying to establish the negative.

The lack of an answer to Open Question~1 is currently a major stumbling block in the classification program for Leavitt path algebras.  The fact we do not have an answer, even in special cases or for elementary examples, indicates there is something important about the structure of simple Leavitt path algebras that we do not yet understand.  In addition, Open Question~1 is not only an impediment for classification of simple Leavitt path algebras, but until we have a solution it essentially impossible to classify nonsimple Leavitt path algebras of finite graphs --- to do so, we would most likely need to deal with the simple ideals and quotients, which are as of yet unmanageable.  Consequently, a solution to Open Question~1 is of paramount importance for the classification program for Leavitt path algebras.

Open Question~1 is compelling to the mathematical community not only for its applications to classification of algebras, but also because whatever the answer turns out to be, it will have consequences for the subjects of Algebra, Functional Analysis, and Dynamics.  As the author sees it, there are three possible answers to Open Question 1: ``Yes", ``No", and ``Sometimes".

If the answer is ``Yes", then this would provide further compelling evidence for the existence of some sort of ``Rosetta Stone" allowing for the translation of results between graph $C^*$-algebras and Leavitt path algebras. Identifying the reason for these similarities could lead to a deeper understanding of the relationships between $C^*$-algebras and algebras.  It could even serve as a call to action for more collaboration between algebraists and analysts.  Perhaps we can find conditions under which dense $*$-subalgebras of $C^*$-algebras have structural properties similar to their ambient $C^*$-algebras.  Perhaps one can find larger classes of algebras for which analogues of $C^*$-algebra results can be proven.  Or perhaps (if we dream big) a version of the Kirchberg-Phillips classification theorem could be proved for a large class of purely infinite simple algebras.

If the answer is ``No", meaning the sign of the determinant is an invariant of Morita equivalence, then the pair $(K_0^\textnormal{alg} (L_K(E)), \sgn (\det (1-A^t_F)))$ would be a complete Morita equivalence invariant for simple Leavitt path algebras of finite graphs.  This would imply that for finite strongly connected graphs, the Morita equivalence class of the Leavitt path algebra coincides exactly with the flow equivalence class of the graph's shift space (cf.~Theorem~\ref{flow-equivalence-thm}).  Consequently, we would have that the Leavitt path algebras and shift spaces are intimately related, suggesting that there is some deeper, not yet understood connection between the algebras and the flow dynamics.

If the answer is ``Sometimes", meaning that for certain graphs changing the sign of the determinant (or performing a Cuntz splice) changes the Morita equivalence class of the  Leavitt path algebra but for other graphs it does not, then we will need to identify exactly which graphs are affected.  This would be the most surprising (and hence for a mathematician the most interesting!) outcome to this question.  If indeed the answer does turn out to be ``Sometimes", this outcome will likely motivate the creation and development of new tools and require the collaboration of algebraists, analysts, and dynamicists to investigate the phenomena that occur.

\section{Classification of Leavitt Path Algebras of Infinite Graphs}
\label{LPA-inf-sec}

As we saw in the previous section, the sign of the determinant condition (and unknown effect of the Cuntz splice) creates an impediment to classifying simple Leavitt path algebras of finite graphs

However, if we continue to look to graph $C^*$-algebras for inspiration, we see that S\o rensen's classification of unital $C^*$-algebras of infinite graphs avoided the sign of the determinant and no Cuntz splice move was needed.  One could therefore hope for a similar classification, using S\o rensen's techniques, for unital Leavitt path algebras of infinite graphs.  (Such graphs have a finite number of vertices and an infinite number of edges.)  This collection of graphs, while avoiding the Cuntz splice, introduces a new problem: What is our candidate for the complete Morita equivalence invariant?  For $C^*$-algebras of finite graphs we used the $K_0$-group, and when we considered Leavitt path algebras of finite graphs, we were in the fortunate situation that for any graph $E$ we have $K_0^\textnormal{alg} (L_K(E)) \cong K_0(C^*(E))$.  However, for graphs with infinitely many edges S\o rensen now had to include the $K_1$-group of the $C^*$-algebra.  In general, for a graph $E$ one has that $K_1^\textnormal{alg} (L_K(E))$ and $K_1(C^*(E))$ are not equal.  In addition, due to Bott periodicity, a $C^*$-algebras really only has only two $K$-groups: the $K_0$-group and then $K_1$-group.  This means that by using $K_0$ and $K_1$, S\o rensen was including \emph{all} the $K$-groups of the graph $C^*$-algebra in the invariant.  For Leavitt path algebras there is no periodicity and the algebraic $K$-groups $K_n^\textnormal{alg} (L_K(E))$ may all be distinct.  Furthermore, for $n \geq 1$ one has that $K_n^\textnormal{alg} (L_K(E))$ depends on the underlying field of the Leavitt path algebra.  This raises the question as to which of the algebraic $K$-groups $K_n^\textnormal{alg} (L_K(E))$ should be included in the invariant.  $K_0^\textnormal{alg}$ and $K_1^\textnormal{alg}$ only?  Some finite number of algebraic $K$-groups? All algebraic $K$-groups?  Should the number of $K$-groups included depend on the field?

Inspired by the techniques of S\o rensen in \cite{Sor}, Ruiz and the author looked for an invariant that would provide the moves needed between the graphs without worrying about whether this invariant involved the algebraic $K$-groups.  It was found that a complete Morita equivalence invariant is provided by the pair $(K_0^\textnormal{alg} (L_K(E)), |E^0_\textnormal{sing}|)$.  Here $|E^0_\textnormal{sing}|$ is the cardinality of the set of singular vertices $E^0_\textnormal{sing}$.  (Recall that a vertex is singular if it either emits no edges or an infinite number of edges.)  Ruiz and the author proved the following in \cite{RT}.

\begin{theorem}[Ruiz and Tomforde] \label{gen-classification-LPA-thm}
Let $K$ be a field, and let $E$ and $F$ be graphs with a finite number of vertices and an infinite number of edges with the property that $L_K(E)$ and $L_K(F)$ are simple.  Then the following are equivalent:
\begin{itemize}
\item[(1)] \ $L_K(E)$ is Morita equivalent to $L_K(F)$.
\item[(2)] \ $K_0^\textnormal{alg}(L_K(E)) \cong K_0^\textnormal{alg}(L_K(F))$ and $|E^0_\textnormal{sing}| = |F^0_\textnormal{sing}|$.
\item[(3)] \ The graph $E$ may turned into the graph $F$ by finitely many applications of Moves (S), (O), (I), (R), their inverses.
\end{itemize}
\end{theorem}

While this result shows that $(K_0^\textnormal{alg} (L_K(E)), |E^0_\textnormal{sing}|)$ is a complete Morita equivalence invariant for unital simple Leavitt path algebras of infinite graphs, this invariant is unsatisfying because it depends on the choice of the graph used to represent the Leavitt path algebra, rather than only on intrinsic properties of the algebra itself.  We prefer to have an invariant based solely on algebraic properties, and if our goal is to lay groundwork for classification of larger classes of algebras, we hope that we can obtain an invariant described entirely in terms of $K$-theory.

Thus we ask whether or not some collection of the algebraic $K$-groups can provide a complete Morita equivalence invariant for unital simple Leavitt path algebras of infinite graphs --- or, equivalently, whether $|E^0_\textnormal{sing}|$ can be determined from some collection of the algebraic $K$-groups.

In \cite{RT} Ruiz and the author showed that in certain situations the answer is ``Yes", and surprisingly the answer depends on the underlying field.  (This is explained by the fact that the higher algebraic $K$-groups of the Leavitt path algebra depend significantly on the underlying field.)  To describe the manageable fields, we need to introduce a bit of terminology.

\begin{definition} \label{no-free-def}
if $G$ is an abelian group, we say $G$ \emph{has no free quotients} if no nonzero quotient of $G$ is a free abelian group.  If $K$ is a field, we say $K$ \emph{has no free quotients} if the multiplicative abelian group $K^\times := K \setminus \{ 0 \}$ has no free quotients.
\end{definition}

\noindent It is shown in \cite[Proposition~6.10]{RT} that the following are all examples of fields with no free quotients.
\begin{itemize}
\item All fields $K$ such that $K^\times$ is a torsion group.
\item All fields $K$ such that $K^\times$ is weakly divisible. 
\item All algebraically closed fields.
\item All fields that are perfect with characteristic $p >0$.
\item All finite fields.
\item The field $\mathbb{C}$ of complex numbers.
\item The field $\mathbb{R}$ of real numbers.
\end{itemize}

\noindent The field $\mathbb{Q}$ is not a field with no free quotients, because $\mathbb{Q}^\times \cong \mathbb{Z}_2 \oplus \mathbb{Z} \oplus \mathbb{Z} \oplus \ldots$.

In \cite{RT} Ruiz and the author showed that when the underlying field has no free quotients, the pair 
$(K_0^\textnormal{alg} (L_K(E)), K_1^\textnormal{alg} (L_K(E)))$ is a complete Morita equivalence invariant.  We emphasize that this includes the case when the underlying field is the complex numbers.

\begin{theorem}[Ruiz and Tomforde] \label{NFQ-LPA-thm}
Let $K$ be a field with no free quotients (see Definition~\ref{no-free-def}), and let $E$ and $F$ be graphs with a finite number of vertices and an infinite number of edges, and with the property that $L_K(E)$ and $L_K(F)$ are simple.  Then the following are equivalent:
\begin{itemize}
\item[(1)] \ $L_K(E)$ is Morita equivalent to $L_K(F)$.
\item[(2)] \ $K_0^\textnormal{alg}(L_K(E)) \cong K_0^\textnormal{alg}(L_K(F))$ and $K_1^\textnormal{alg}(L_K(E)) \cong K_1^\textnormal{alg}(L_K(F))$.
\item[(3)] \ The graph $E$ may turned into the graph $F$ by finitely many applications of Moves (S), (O), (I), (R), their inverses.
\end{itemize}
\end{theorem}

Moreover, in \cite[\S11]{RT} Ruiz and the author produce an example of graphs $E$ and $F$ with finitely many vertices, infinitely many edges, and having the following properties: $L_\mathbb{Q}(E)$ and $L_\mathbb{Q}(F)$ are simple, $K_0^\textnormal{alg}(L_\mathbb{Q}(E)) \cong K_0^\textnormal{alg}(L_\mathbb{Q}(F))$, $K_1^\textnormal{alg}(L_\mathbb{Q}(E)) \cong K_1^\textnormal{alg}(L_\mathbb{Q}(F))$, and $K_2^\textnormal{alg}(L_\mathbb{Q}(E)) \not\cong K_2^\textnormal{alg}(L_\mathbb{Q}(F))$.  Thus the $K_0^\textnormal{alg}$-groups and $K_1^\textnormal{alg}$-groups of $L_\mathbb{Q}(E)$ and $L_\mathbb{Q}(F)$ are isomorphic, but $L_\mathbb{Q}(E)$ and $L_\mathbb{Q}(F)$ are not Morita equivalent.  Hence the pair of the $K_0^\textnormal{alg}$-group and $K_1^\textnormal{alg}$-group can fail to be a complete Morita equivalence invariant when the underlying field is not a field with no free quotients.

This raises the question of whether higher algebraic $K$-groups can be included to produce a complete Morita equivalence invariant for other fields.  This was answered affirmatively for number fields in \cite{GRTW}.  (Recall that a number field is a field that is a finite extension of $\mathbb{Q}$.)  The following was proven in \cite{GRTW}.

\begin{theorem}[Gabe, Ruiz, Tomforde, and Whalen] \label{number-fields-LPA-thm}
Let $K$ be a number field, and let $E$ and $F$ be graphs with a finite number of vertices and an infinite number of edges, and with the property that $L_K(E)$ and $L_K(F)$ are simple.  Then the following are equivalent:
\begin{itemize}
\item[(1)] \ $L_K(E)$ is Morita equivalent to $L_K(F)$.
\item[(2)] \ $K_0^\textnormal{alg}(L_K(E)) \cong K_0^\textnormal{alg}(L_K(F))$ and $K_6^\textnormal{alg}(L_K(E)) \cong K_6^\textnormal{alg}(L_K(F))$.
\item[(3)] \ The graph $E$ may turned into the graph $F$ by finitely many applications of Moves (S), (O), (I), (R), their inverses.
\end{itemize}
\end{theorem}

This shows that the pair $(K_0^\textnormal{alg}(L_K(E)), K_6^\textnormal{alg}(L_K(E)))$ is a complete Morita equivalence invariant for these Leavitt path algebras when the field is a number field.  In addition, since $\mathbb{Q}$ is a number field, this result covers the example produced by Ruiz and the author in \cite[\S11]{RT}.

The situation for other fields is unclear, and this leads us to an important open question.

$ $

\noindent \textbf{Open Question 2:} \emph{Let $K$ be a field, and let $E$ and $F$ be graphs with a finite number of vertices and an infinite number of edges, and with the property that $L_K(E)$ and $L_K(F)$ are simple. If $K_n^\textnormal{alg} (L_K(E)) \cong K_n^\textnormal{alg} (L_K(F))$ for all $n \in \mathbb{N} \cup \{ 0 \}$, then is it the case that $L_K(E)$ is Morita equivalent to $L_K(F)$?}

$ $

\noindent In light of Theorem~\ref{gen-classification-LPA-thm}, Open Question~2 is equivalent to asking the following: ``If $K$ is a field and $E$ is a graph with a finite number of vertices and an infinite number of edges, and for which $L_K(E)$ is simple, then is it possible to determine $|E^0_\textnormal{sing}|$ from the set of algebraic $K$-groups $\{ K_n^\textnormal{alg} (L_K(E)) : n =0,1,2 \ldots \}$?"

Theorem~\ref{NFQ-LPA-thm} and Theorem~\ref{number-fields-LPA-thm} show that Open Question~2 has an affirmative answer when the field either has no free quotients or is a number field.  This means that for simple Leavitt path algebras over these fields the only missing part of a classification is to answer Open Question~1 and determine whether the sign of the determinant is a Morita equivalence invariant in the case of finite graphs.

To obtain a classification of all simple Leavitt path algebras by algebraic $K$-theory, a positive answer to Open Question~2 is necessary.  If the answer to Question~2 is negative in general, then a general classification in terms of algebraic $K$-theory will not be possible and we will need to restrict our attention to Leavitt path algebras over particular fields.  Fortunately, for many fields of interest (e.g., $\mathbb{C}$, $\mathbb{R}$, finite fields, $\mathbb{Q}$, number fields) we already know that Open Question~2 has a positive answer.

\begin{acknowledgement}
The author thanks the organizers of the 2015 Abel Symposium, Christian Skau (Norwegian University of Science and Technology), Toke M. Carlsen (University of the Faroe Islands), Nadia Larsen (University of Oslo) and Sergey Neshveyev (University of Oslo) for their hospitality and the opportunity to attend.  This work, including the author's travel to the Abel Symposium, was supported by a grant from the Simons Foundation (\#210035 to Mark Tomforde).
\end{acknowledgement}


\begin{thebibliography}{99.}%

\bibitem{AAP}
G.~Abrams and G.~Aranda Pino, \emph{The Leavitt path algebra of a graph}, J. Algebra \textbf{293} (2005), no. 2, 319--334.

\bibitem{ALPS}  
G. Abrams, A. Louly, E. Pardo, and C. Smith, \emph{Flow invariants in the classification of Leavitt path algebras}, J.~Algebra \textbf{333} (2011),  202--231.

\bibitem{Cun2}
J.~Cuntz, \emph{$K$-theory for certain $C^*$-algebras},
Ann. Math. \textbf{113} (1981), 181--197.

\bibitem{CK}
J.~Cuntz and W.~Krieger, \emph{A class of $C^*$-algebras and
topological Markov chains}, Invent. Math. \textbf{56} (1980),
251--268.

\bibitem{Dri}
D.~Drinen, \emph{Viewing AF-algebras as graph algebras},
Proc. Amer. Math. Soc. \textbf{128} (2000), 1991--2000. 

\bibitem{ERR}
S.~Eilers, G.~Restorff, and E.~Ruiz, \emph{Classification of graph $C^*$-algebras with no more than four primitive ideals}, Operator algebra and dynamics, 89--129, Springer Proc. Math. Stat., 58, Springer, Heidelberg, 2013.

\bibitem{ERRS}
S.~Eilers, G.~Restorff, E.~Ruiz, and A.~P.~W.~S\o rensen, \emph{Geometric classification of unital graph $C^*$-algebras of real rank zero}, preprint, 2015. \texttt{arXiv:1505.06773 [math.OA]}

\bibitem{ET1}
S.~Eilers and M.~Tomforde, \emph{On the classification of nonsimple graph $C^*$-algebras}, Math. Ann. \textbf{346} (2010), 393--418. 

\bibitem{ET}
G.A.~Elliott, and A.~Toms, \emph{Regularity properties in the classification program for separable amenable $C^*$-algebras}, Bull. Amer. Math. Soc. (N.S.) \textbf{45} (2008), no. 2, 229--245.

\bibitem{GRTW}
J.~Gabe, E.~Ruiz, M.~Tomforde, and T.~Whalen, \emph{$K$-theory for Leavitt path algebras: computation and classification}, J. Algebra \textbf{433} (2015), 35--72.

\bibitem{JS}
R.~Johansen and A.~P.~W.~S\o rensen, \emph{The Cuntz splice does not preserve $*$-isomorphism of Leavitt path algebras over $\mathbb{Z}$}, preprint, 2015.  \texttt{arXiv:1507.01247 [math.RA]}

\bibitem{LM}
D.~Lind and B.~Marcus, An Introduction to Symbolic
Dynamics and Coding, Cambridge University Press, Cambridge, 1995.

\bibitem{MN}
R.~Meyer and R.~Nest, Ryszard, \emph{$C^*$-algebras over topological spaces: filtrated K-theory}, Canad. J. Math. \textbf{64} (2012), no. 2, 368--408.

\bibitem{Phi}
N.~C.~Phillips, \emph{A classification theorem for nuclear
purely infinite simple $C^*$-algebras},  Doc. Math.
\textbf{5} (2000), 49--114. 

\bibitem{Ro}
M.~R\o rdam, \emph{Classification of Cuntz-Krieger algebras},
$K$-theory \textbf{9} (1995), 31--58. 

\bibitem{Ror-book}
M.~R\o rdam, Classification of Nuclear, Simple $C^*$-algebras,
Encyclopaedia of Mathematical Sciences, vol. 126, Springer, Berlin,
2001.

\bibitem{Ror}
M.~R\o rdam, \emph{Structure and classification of $C^*$-algebras}, International Congress of Mathematicians. Vol. II, 1581--1598, Eur. Math. Soc., Z\"urich, 2006. 

\bibitem{RT}
E.~Ruiz, and M.~Tomforde, \emph{Classification of unital simple Leavitt path algebras of infinite graphs}, J. Algebra \textbf{384} (2013), 45--83.

\bibitem{Sor}
A.~P.~W.~S\o rensen, \emph{Geometric classification of simple graph algebras}, Ergodic Theory Dynam. Systems \textbf{33} (2013), no. 4, 1199--1220.

\bibitem{Szy}
W.~Szyma\'nski, \emph{The range of $K$-invariants for 
$C^*$-algebras of infinite graphs}, Indiana Univ. Math. J.
\textbf{51} (2002), 239--249.

\bibitem{TWW}
A.~Tikuisis, S.~White, and W.~Winter,
\emph{Quasidiagonality of nuclear $C^*$-algebras}, preprint (2015).

\bibitem{Tom}
M.~Tomforde, \emph{Leavitt path algebras with coefficients in a commutative ring}, J. Pure Appl. Algebra
\textbf{215} (2011), 471--484. 




\end{thebibliography}
%

\end{document}